\documentclass{article}

\usepackage{amsmath, amssymb,amsthm, amstext}
\usepackage{mathrsfs}
\usepackage{graphicx}
\usepackage{dsfont}
\usepackage{slashbox}

\newtheoremstyle{break}
   {9pt}
   {9pt}
   {\itshape}
   {}
   {\bfseries}
   {.}
   {0.5em}
   {}

\newtheoremstyle{definitionbreak}
   {9pt}
   {9pt}
   {\rm}
   {}
   {\bfseries}
   {.}
   {0.5em}
   {}

\theoremstyle{break}
\newtheorem{thm}{Proposition}[section]
\newtheorem{them}[thm]{Theorem}
\newtheorem{lem}[thm]{Lemma}
\newtheorem{cor}[thm]{Corollary}

\theoremstyle{definitionbreak}
\newtheorem{definition}[thm]{Definition}
\newtheorem{rem}[thm]{Remark}
\newtheorem{notat}[thm]{Notation}
\newtheorem{r&n}[thm]{Remark and Notation}
\newtheorem{ex}[thm]{Example}
\newtheorem{algor}[thm]{Algorithm}

\def\q{{\mathfrak{q}}}
\def\p{{\mathfrak{p}}}
\def\a{{\mathfrak{a}}}

\def\L{{\mathfrak L}}
\def\E{{\mathfrak E}}
\def\lra{\longrightarrow}

\def\Epi{\twoheadrightarrow}
\def\Llra{\Longleftrightarrow}

\def\inkl{\hookrightarrow}
\def\bs{\backslash}
\def\ss{\subseteq}

\def\nss{\not\subseteq}
\def\Nat{\mathbb{N}_0}
\def\ZZ{\mathbb{Z}}
\def\NN{\mathbb{N}}
\def\PP{\mathbb{P}}

\def\LL{\mathbb{L}}
\def\EE{\mathbb{E}}
\def\defin{\mathrel{\mathop:}=}
\def\ON{\setminus\{0\}}

\def\>>{>\!\!>}
\def\<<{<\!\!<}
\def\sat{^{\mbox{{\footnotesize\rm sat}}}}

\DeclareMathOperator{\Ass}{Ass}
\DeclareMathOperator{\Var}{Var}

\DeclareMathOperator{\NZD}{NZD}

\DeclareMathOperator{\hth}{heigth}
\DeclareMathOperator{\im}{im}

\DeclareMathOperator{\length}{length}
\DeclareMathOperator{\LC}{LC}
\DeclareMathOperator{\LT}{LT}

\DeclareMathOperator{\Proj}{Proj}
\DeclareMathOperator{\mProj}{mProj}

\DeclareMathOperator{\I}{I}

\DeclareMathOperator{\In}{In}
\DeclareMathOperator{\Sec}{Sec}
\DeclareMathOperator{\Jn}{Join}

\newcommand{\res}{\!\!\upharpoonright}
\newcommand{\ol}[1]{\overline{#1}}

\newcommand{\K}{{\cal K}}

\begin{document}

\title{Partial elimination ideals and secant cones}
\author{Simon Kurmann}
\date{}

\maketitle

\begin{abstract}
For any $k \in \Nat$, we show that the cone of $(k+1)$-secant lines of a closed subscheme $Z \ss \PP^n_K$ over an algebraically closed field $K$ running through a closed point $p \in \PP^n_K$ is defined by the $k$-th partial elimination ideal of $Z$ with respect to $p$. We use this fact to give an algorithm for computing secant cones. Also, we show that under certain conditions partial elimination ideals describe the length of the fibres of a multiple projection in a way similar to the way they do for simple projections. Finally, we study some examples illustrating these results, computed by means of {\sc Singular}.
\end{abstract}

\section{Introduction}	

Partial elimination ideals (PEIs) have been introduced by M. Green in \cite{G} in relation to generic initial ideals. In this article we give a definition of partial elimination ideals that is independent of the choice of coordinates and can be used to study simple projections with arbitrary centre as well as multiple projections with certain nice properties. Finally, we give an algorithm which utilises PEIs to compute the secant cone of a projective scheme with respect to an arbitrary point.\\
Let $K$ be an algebraically closed field, let  $n \in \Nat$, and let $R \defin K[x_0, \ldots, x_n]$ be the polynomial ring in $n+1$ indeterminates. Let $\PP^n = \PP^n_K \defin \Proj(R)$; let $\mProj(R)$ denote the set of maximal ideals in $\Proj(R)$, that is the set of closed points of $\PP^n$. Let $\p \in \mProj(R)$ be a closed point, let $S \defin K[\p_1]$ be the homogeneous $K$-subalgebra of $R$ generated by the linear forms of $\p$, and let $\a \ss R$ be a graded ideal. Then, for $k \in \Nat$ and each $y \in R_1\bs\p_1$, the $k$-th partial elimination ideal of $\a$ with respect to $\p$ is \[\K^\p_k(\a) = \bigoplus_{d \in \ZZ} \{f \in S_d \mid \exists g \in (\p^{d+1})_{d+k}: y^kf+g \in \a_{d+k}\},\] which is a graded ideal in $S$. Note that $\K^\p_k(\a)$ is independent of our choice of $y \in R_1\bs\p_1$ (Corollary \ref{JustOneCoord}). This definition is indeed a generalisation of that in \cite{G}, where PEIs are defined in the case where $\p =  (x_1,\ldots,x_n)$ and $y = x_0$. In \cite{C-S} it is shown that in this situation a Gr{\"o}bner basis of $\K^{\p}_k(\a)$ is given by the leading coefficients in $x_0$ of those elements of a Gr{\"o}bner basis of $\a$ whose degree in $x_0$ is less or equal than $k$ (with respect to an elimination ordering on $R$). As taking PEIs commutes with coordinate transformation (Lemma \ref{CoordTransf}) it is therefore quite easy to compute PEIs.\\
Now, let $Z \ss \PP^n$ be the closed subscheme defined by a homogeneous ideal $I_Z \ss R$ such that $p = (1:0:\cdots:0) \notin Z$, and let $\pi:\PP^n\bs\{p_0\} \to \PP^{n-1}$ be the projection from $p$ to the subspace $\PP^{n-1} \ss \PP^n$ whose homogeneous coordinate ring $S = K[x_1, \ldots, x_n]$ is generated by the linear forms of the ideal of $p$. Then it is well known (see \cite[Proposition 6.2]{G}, \cite[Theorem 3.5]{C-S}) that $\K^{\p}_k(I_Z)$ is the vanishing ideal of the set $\{q \in \pi(Z) \mid l(Z \cap \langle q,p\rangle) > k\}$, where $\langle q,p\rangle \ss \PP^n$ is the line spanned by $q$ and $p$ and $l(Z \cap \langle q,p\rangle)$ denotes the length of the fibre $(\pi\res_Z)^{-1}(q) = Z\cap\langle q,p\rangle$ over $q$. (Incidentally, this is the reason why we consider schemes instead of varieties, as we have  to study intersection length and closed points counted with multiplicity.) Theorem \ref{Theorem} gives a sligthly more general version of this result, and we give a proof that demonstrates the relation between PEIs and homogeneous elements of the homogeneous ring $R/I_Z$ which behave analogously to superficial elements in local rings. Further, Theorem \ref{Theorem} immediately gives rise to Proposition \ref{SecPEI} which states that $\sqrt{\K^\p_k(I_Z)R}$ is the (scheme theoretic) ideal of the $(k+1)$-secant cone $\Sec^{k+1}_p(Z)$ of $Z$ with respect to the closed point $p \in \PP^n\bs Z$ corresponding to any $\p \in \mProj(R)$. This allows us to define an algorithm for computing secant cones and secant loci of a projective subscheme $Z \ss \PP^n$ with respect to a closed point $p\in \PP^n\bs Z$ (see Algorithm \ref{SecComp}):
\medskip\\
{\bf Input:} Ideal $I_Z \ss R$ of $Z$, ideal $p \in \mProj(R)$.\\
1. Define a linear coordinate transformation $\psi: R \stackrel{\cong}{\lra} R$ such that $\psi(p) = (x_1, \ldots,x_n)$.\\
2. Compute a Gr{\"o}bner basis $G$ of $\psi(I_Z)$ with respect to an elimination ordering on $R$.\\
3. Set $k_0 \defin \max\{\deg_{x_0}(g) \mid g \in G\}$.\\
4. For $0 \leq k \leq k_0$, set $G_k \defin \{\LT_{x_0}(g) \mid g \in G \wedge \deg_{x_0}(g) < k\}$.\\
5. Compute $\K^p_k(I_Z)R = \psi^{-1}(G_k)R$ for $0 \leq k \leq k_0$.\\
6. Compute $\sqrt{\K^p_k(I_Z)R}$ for $0 \leq k \leq k_0$.\\
{\bf Output:} Ideals $\sqrt{\K^p_0(I_Z)R}, \ldots, \sqrt{\K^p_{k_0}(I_Z)R}$ of $\Sec_p^1(Z), \ldots, \Sec^{k_0+1}_p(Z)$.
\medskip\\
In section 2 we give a sligthly different definition of partial elimination ideals and show that this definition indeed describes what we were looking for and that it is independent of the choice of coordinates. In section 3 we formulate the main results about PEIs with respect to simple projections and secant cones. In section 4 we give two results about PEIs and multiple projections useful for the consideration of examples. Section 5 explains the algorithm for computing secant cones and secant loci using PEIs. Finally, in section 6 we consider some examples.
\medskip\\
{\bf Acknowledgment.} I would like to thank my Ph.D. advisor Markus P. Brodmann for suggesting that I write this paper; without his support it would not have been written. I thank Kangjing Han for some stimulating discussions about the geometry of PEIs. Peter Schenzel had some kind words and comments about my first effort in programming which I found quite motivating, and Euisung Park provided an easy answer to a geometric problem of mine. Finally, I thank Felix Fontein, Thomas Preu, Fred Rohrer, and Rodney Y. Sharp for proofreading and helpful suggestions.

\section{Partial Elimination Ideals}\label{DefPEI}

\begin{notat} (A) For the remainder of this article, by $p \in \PP^n$ we always mean a closed point $p$ of the projective $n$-space $\PP^n_K$. Moreover, we identify $\PP^n$ and $\mProj(R)$ and just write $p \in \mProj(R)$ for the homogeneous ideal of $p \in \PP^n$. If we want to emphasize an algebraic point of view, we use gothic letters such as $\q$ and $\p$, which always relate to closed points denoted by the according latin letters such as $q$ and $p$.\\
(B) Let $\p \in \mProj(R)$ and $K[\p_1] \ss R$ be the graded $K$-subalgebra of $R$ generated by the linear forms of $\p$. Frequently, we just write $S \defin K[\p_1]$ and consider $S$ as the homogeneous ring of some linear subspace $\PP^{n-1} = \PP^{n-1}_K(\p) \defin \Proj(K[\p_1]) \ss
\PP^n$. Now, let $y \in R_1\bs\p_1$. Then, the $K$-vectorspace $R_1$ is generated by $y$ and $\p_1$ and
therefore $R = K[y,\p_1] = S[y]$. Let $f \in R$. We may consider $f$ as a polynomial in $y$ over $S$ and write $\deg_{y}(f)$ for the degree of $f \in S[y]$ in
$y$. Furthermore, let $\LC_{y}(f)$ and $\LT_{y}(f)$ denote the leading coefficient and the leading term of $f$, respectively, as a polynomial in $y$ over $S$. Note that using this notations we are not considering the standard grading on $R$ but the one induced by $R = S[y]$.
\end{notat}

\begin{definition}\label{Def} Let $\a \ss R$ be a graded ideal, let $\p \in \mProj(R)$, and let $k \in \Nat$. We define the {\it $k$-th partial elimination ideal} (abbreviated PEI) {\it of $\a$ with respect to $\p$} by \[\K^\p_k(\a) \defin \{f \in K[\p_1] \mid \forall y \in R_1\bs\p_1 \, \exists g \in R: \deg_y(g)<k \wedge y^kf+g \in \a\}.\] $\K^\p_k(\a)$ is a graded ideal of $S$ whose $d$-th graded component is given by \[\setlength\arraycolsep{0.2em} \begin{array}{rcl} \K_k^\p(\a)_d &=& \{f \in S_d \mid \forall y \in R_1\bs\p_1 \, \exists g \in R_{d+k}: \deg_y(g)<k \wedge y^kf+g \in \a_{d+k}\}\\ &=& \{f \in S_d \mid \forall y \in R_1\bs\p_1 \, \exists g \in (\p^{d+1})_{d+k}: y^kf+g \in \a_{d+k}\cap (\p^d)_{d+k}\}.\end{array}\] Finally, we define $\K_{-1}^\p(\a) = 0$. In this way we get an ascending chain of graded ideals of $S$ \[0 = \K_{-1}^\p(\a) \ss \a \cap S = \K_0^\p(\a) \ss \K_1^\p(\a) \ss \cdots \ss \K^\p_k(\a) \ss \K^\p_{k+1} \ss \cdots\]
\end{definition}

\begin{notat}
For the remainder of this section, we fix a graded ideal $\a \ss R$ and a closed point $\p \in \mProj(R)$. For all $k \in \ZZ$, let \[\widetilde{\K}^\p_k(\a) \defin
\bigoplus_{d\in\ZZ}\a_d\cap(\p^{d-k})_d,\] where $\p^d = R$ for $d\leq0$. $\widetilde{\K}^\p_k(\a)$ is a graded $S$-module, and it holds $\widetilde{\K}_{k-1}^{\p}(\a) \ss \widetilde{\K}_k^{\p}(\a)$ for all $k \in \ZZ$. Moreover, a straightforward calculation shows that \[\forall y \in R_1\bs\p_1: \widetilde{\K}_k^\p(\a) = \{f \in \a \mid \deg_{y}(f)\leq k\},\] meaning that for each $y \in R_1\bs\p_1$ we can write any element $f \in \widetilde{\K}^\p_k(\a)$ uniquely as $f = y^kf_0 + g$ with $f_0 \in S$ and $g \in R$ such that $\deg_{y}(g) < k$. 
\end{notat}

\begin{lem}\label{Isom} For all $k \in \Nat$ and all $y \in R_1\bs\p_1$, there is an isomorphism of graded $S$-modules \[\varphi^{y}_k(\a):
\widetilde{\K}_k^\p(\a)/\widetilde{\K}_{k-1}^\p(\a)(-k) \stackrel{\cong}{\lra} \K_k^\p(\a), \ol{f} = y^kf_0 +g + \widetilde{\K}_{k-1}^\p(\a) \mapsto f_0.\]
\end{lem}
	
\begin{proof} Let $y \in R_1\bs\p_1$, and let $k \in \Nat$. There is a morphism of graded $S$-modules \[\ol{\varphi}^{y}_k(\a): \widetilde{\K}_k^\p(-k) \lra S, f = y^kf_0+g \mapsto f_0.\] By definition we find $\K_k^\p(\a) \ss \im(\ol{\varphi}^{y}_k(\a))$; we want to show that the reverse inclusion holds, too. So, let $y' \in R_1\bs\p_1$ be arbitrary, let $d \in \Nat$, and let $f_0 \in \im(\ol{\varphi}^{y}_k(\a))_{d} \ss S_d$. Then $f_0 \in (\p^d)_d \ss R_d$, and there is an element $g \in (\p^{d+1})_{d+k}$ such that $y^kf+g \in \a_{d+k} \cap (\p^{d})_{d+k}$. As $R_1$ is generated by $y$ and $\p_1$ over $K$, we find elements $\lambda \in K\ON$ and $v \in \p_1$ such that $y = \lambda y' + v$, that is $y^kf_0 = \lambda^ky'^kf_0 + uf_0$ for some $u \in \p$. Therefore, $\frac{uf_0}{\lambda^k}+\frac{g}{\lambda^k} \in \p^{d+1}$ and $y'^kf_0 + \frac{uf_0}{\lambda^k}+\frac{g}{\lambda^k} \in \a$, proving indeed $f_0 \in \K_k(\a)$.
On the other hand, it is easy to see that $\ker(\ol{\varphi}^{y}_k(\a)) = \widetilde{\K}_{k-1}^\p(\a)(-k)$. This immediately gives the desired isomorphism.
\end{proof}

\begin{cor}\label{JustOneCoord} For all $y \in R_1\bs\p_1$ and all $k \in \Nat$ it holds \[\setlength\arraycolsep{0.2em}\begin{array}{rcl}\K^\p_k(\a) &=& \bigoplus_{d \in \ZZ}\{f \in S_d \mid \exists g \in R_{d+k}: \deg_{y}(g) < k \wedge y^kf+g \in \a_{d+k}\}\\ &=& \bigoplus_{d \in \ZZ}\{f \in S_d \mid \exists g \in (\p^{d+1})_{d+k}: y^kf+g \in \a_{d+k}\cap (\p^d)_{d+k}\}.\end{array}\]
\end{cor}

\begin{proof} Let $y \in R_1\bs\p_1$, $k \in \Nat$, and write $\K^{y}_k(\a) \defin \bigoplus_{d \in \ZZ}\{f \in S_d \mid \exists g \in (\p^{d+1})_{d+k}: y^kf+g \in \a_{d+k}\cap \p^d_{d+k}\}$. As in the above proof, there is an isomorphism of graded $S$-modules \[\widetilde{\K}_k^\p(\a)/\widetilde{\K}_{k-1}^\p(\a)(-k) \stackrel{\cong}{\lra} \K^{y}_k(\a), y^kf_0 +g + \widetilde{\K}_{k-1}^\p(\a) \mapsto f_0.\] This gives $\K^\p_k(\a) \cong \K^{y}_k(\a)$. As by definition $\K^\p_k(\a) \ss \K^{y}_k(\a)$, we get our claim.
\end{proof}

Note that the formula of Corollary \ref{JustOneCoord} is indeed the same as the one given in the introduction as $S_d \ss \p^d$.

\begin{lem}\label{CoordTransf} Let $\psi: R \stackrel{\cong}{\lra} R$ be a graded ring automorphism, and let $k \in \Nat$. Then \[\psi(\K^\p_k(\a)) = \K^{\psi(\p)}_k(\psi(\a)).\]
\end{lem}

\begin{proof} Let $y \in R_1\bs\p_1$. Then $\psi(y) \in R_1\bs\psi(\p)_1$ and $\langle y,\p_1\rangle_K = R_1 = \psi(R_1) = \langle \psi(y),\psi(\p)\rangle_K$, where $\langle y,\p_1\rangle_K$ denotes the $K$-vectorspace generated by $y$ and $\p_1$. Now let $d \in \ZZ$, and let $f \in S_d$. Then by Corollary \ref{JustOneCoord} we see \[
\setlength\arraycolsep{0.2em}\begin{array}{rclcl} f \in \K^\p_k(\a)_d \ss K[\p_1]_d &\Llra& \exists g \in (\p^{d+1})_{d+k}: y^kf+g \in \a_{d+k}\cap (\p^d)_{d+k}\\ &\Llra& \exists g' \in (\psi(\p)^{d+1})_{d+k}:\\ && \psi(y)^k\psi(f) + g' \in \psi(\a)_{d+k}\cap(\psi(\p)^d)_{d+k}\\ &\Llra& \psi(f) \in \K^{\psi(\p)}_k(\psi(\a))_d \ss K[\psi(\p)_1]_d.\end{array}\]
\end{proof}

\begin{rem} Corollary \ref{JustOneCoord} means that to compute PEIs it is enough to look at one element $y \in R_1\bs\p_1$, while Lemma \ref{CoordTransf} tells us that computing PEIs commutes with coordinate transformations. Definition \ref{Def} therefore indeed gives a generalisation of the partial elimination ideals defined in \cite[6.1]{G} which is independent of a choice of coordinates of $R$.
\end{rem}

\section{Partial Elimination Ideals and Secant Lines}

\begin{notat}
For the remainder of this section, let $\a \ss R$ be a graded ideal, let $\p \in \mProj(R)$, and let $y_0 \in R_1\bs\p_1$, that is $R_1$ is generated by $y_0$ and $\p_1$ over $K$. For a graded ideal $\q \ss S = K[\p_1]$ let \[\ol{S} \defin S/((\a\cap S)+\q).\] According to the homogeneous normalization lemma, if $\dim(\ol{S}) = 1$, then there is an element $y_1 \in S_1$ such that $K[y_1] \ss \ol{S}$ is a finite integral extension (here and later we identify indeterminates and their residue classes if there is no danger of mistakes). Furthermore, we can write \[\ol{R} \defin R/((\a\cap S)R + \q R) = \ol{S}[y_0].\] The ring extension $K[y_0,y_1] \ss \ol{R}$ is finite and integral, too. Let \[\ol{\p} \defin \p/((\a\cap S)R + \q R),\] and let \[\ol{\a} \defin \a +\q R/((\a\cap S)R + \q R).\] Then $\hth(\ol{\a}) = 1 = \hth(\ol{\a} \cap K[y_0,y_1])$, so $(\ol{\a} \cap K[y_0,y_1])\sat$ is a principal ideal.\\
If $\a \nss \p$, it holds $\sqrt{\a+\p} = R_+$, and therefore there exists an integer $t \in \Nat$ and an element $g \in \p_t$ such that $y_0^t+g \in \a$. Hence $\deg_{y_0}(g) < t$ implies $1_S \in \K_t^\p(\a)$. So, if $\a\cap S+\q \neq S$ there exists an integer \[k_0 \defin \max\{k \in \Nat\cup\{-1\} \mid \K^\p_k(\a) \ss \a\cap S + \q\}.\] 
\end{notat}

\begin{lem}\label{asatGen}
 Assume $\a \nss \p$. Let $\q \ss S$ be a graded ideal such that $\dim(\ol{S}) = 1$, and such that $\ol{\a}\sat \ss \ol{R}$ is a principal ideal. Then any generator $\ol{h}$ of $\ol{\a}\sat$ can be written as \[\ol{h} = h_0y_0^{k_0+1}+\ol{g} \in \ol{R}_{k_0+1}\] with $h_0 \in K\bs\{0\}$ and $\ol{g} \in \ol{\p}_{k_0+1}$.
\end{lem}

\begin{proof}
 Let $y_1 \in \ol{S}_1$ such that $K[y_1] \ss \ol{S}$ is finite and integral. Let $\ol{h} \in \ol{R}$ be a homogeneous generator of $\ol{\a}\sat$, and let $l \defin \deg(\ol{h})$. As $\ol{R} = K[y_0,\ol\p_1]$, we can write \[\ol{h} = h_0y_0^l + \ol{g},\] where $h_0 \in K$ and $\ol{g} \in \ol{\p}$. As $\ol{\a}\sat \nss \ol{\p}$, it follows $h_0 \neq 0$.\\
So, we need only show that $l = k_0+1$. Let $g \in \p_l$ be a representative of $\ol{g}$. Then $h \defin h_0y_0^l +g$ is a representative of $\ol{h}$. For $d \>> 0$ it holds $\ol{h}y_1^d \in \ol{\a}_{d+l} \cap (\ol{\p}^d)_{d+l}$ and therefore $hy_1^d \in \a_{d+l} \cap (\p^d)_{d+l} + ((\a\cap S)R + \q R)_{d+l}$. Thus, there are elements $v \in \a_{d+l}\cap (\p^{d})_{d+l}$ and $w \in ((\a\cap S)R + \q R)_{d+l}$ such that $hy_1^d = v+w$. In particular $\LT_{y_0}(v+w) = h_0y_0^ly_1^d$. As $v \in (\p^d)_{d+l}$ it holds $\deg_{y_0}(v)\leq l$ and therefore $\deg_{y_0}(w) \leq l$. We write $w = y_0^lw_0 + \tilde{w}$, where $w_0 \in S_d$ and $\tilde{w} \in R_{d+l}$ with $\deg_{y_0}(\tilde{w}) < l$. As $y_0^t \notin (\a\cap S)R + \q R$ for all $t \in \Nat$ it follows $w_0 \in \a\cap S + \q$, and as $y_1^t \notin \a\cap S + \q$ for all $t \in \Nat$ we finally get $w_0 \neq h_0y_1^d$. This means that $\LT_{y_0}(v) = y_0^l(h_0y_1^d - w_0)$ and hence $h_0y_1^d - w_0 \in \K_l^\p(\a)$. If $l \leq k_0$, this would imply $y_1^d \in \a\cap S + \q$, a contradiction. It follows $l > k_0$.\\
On the other hand let $k \in \Nat$ with $k < l$, let $d \in \Nat$, and let $f \in \K^\p_k(\a)_d$. We want to show that $f \in \a\cap S + \q$. There is an element $g \in (\p^{d+1})_{d+k}$ such that $y_0^kf + g \in \a_{d+k}\cap (\p^d)_{d+k}$, that is \[y_0^kf + g + ((\a\cap S)R + \q R)_{d+k} \in \ol{\a}\sat_{d+k}\cap(\ol{\p}^d)_{d+k} = 0,\] where the last equality holds because of $\deg(\ol{h})>k$ and $\ol{h} \notin \ol{\p}$. 
So $y_0^kf + ((\a\cap S)R + \q R)_{d+k} = -g + ((\a\cap S)R + \q R)_{d+k}$. If we assume $f \notin ((\a\cap S)R + \q R)_{d+k}$, we therefore immediately get the contradiction \[\setlength\arraycolsep{0.2em}\begin{array}{lclclcl}k &=& \deg_{y_0}(y_0^kf+ ((\a\cap S)R + \q R)_{d+k})\\ &=& \deg_{y_0}(-g + ((\a\cap S)R + \q R)_{d+k}) &\leq& \deg_{y_0}(g) &<& k.\end{array}\] This proves $\K_k^\p(\a) \ss \a\cap S + \q$ for all $k < l$ and therefore $l \leq k_0+1$.
\end{proof}

\begin{r&n} (A) Let $M$ be a finitely generated graded $R$-module, $M \neq 0$, let $p_M$ denote the Hilbert polynomial of $M$, and let $d \defin \dim(M)$. Then, we denote the Hilbert multiplicity of $M$ by \[e_0(M) \defin \left\{\begin{array}{ll}\length(M) & d = 0\\ (d-1)!\cdot\LC(p_M) & d >0\end{array}\right. .\]
(B) Let $M$ be a finitely generated graded $R$-module with $d \defin \dim(M)>0$, and let $r \in \{1,\ldots, d-1\}$. Let $t_1,\ldots, t_r \in \NN$, and let $h_i \in R_{t_i}$ for $1 \leq i \leq r$ such that $\dim(M/\Sigma_{i=1}^rh_iM) = d-r$, that is $h_1, \ldots, h_r$ form a system of homogeneous parameters of $M$. Then \[e_0(M/\Sigma_{i=1}^rh_iM) \geq t_1\cdots t_r\cdot e_0(M).\] Furthermore, the following conditions are equivalent: 
\begin{enumerate}
 \item[(a)] $e_0(M/\Sigma_{i=1}^rh_iM) = t_1\cdots t_r\cdot e_0(M)$;
 \item[(b)] $\forall s \in \{1, \ldots, r\}: h_s \notin \bigcup\{\p \in \Ass_R(M/\Sigma_{i=1}^{s-1}h_iM) \mid \dim(R/\p) \geq d-s\}$.
\end{enumerate}
If these two equivalent conditions hold, we call $h_1, \ldots, h_r$ a {\it system of multiplicity parameters of degree $t_1, \ldots, t_r$ for $M$}. By (b) it follows that every $M$-sequence is a system of multiplicity parameters. If $r = 1$, we just call $h = h_1$ a {\it multplicity parameter of degree $t = t_1$ for $M$.} 
Multiplicity parameters are the analogue in homogeneous rings of superficial elements in local rings (see \cite[VIII \S 7.5]{B}.)
\end{r&n}

\begin{them}\label{Theorem} Assume $\a \nss \p$. Let $\q \ss S$ be a graded ideal such that $\dim(\ol{S}) = 1$ and such that $\ol{\a}\sat \ss \ol{R}$ is a principal ideal. Then \[e_0(R/(\a + \q R)) = (k_0+1)\cdot e_0(S/(\a \cap S + \q)).\]
\end{them}

\begin{proof} Let $y_1 \in \ol{S}_1$ such that $K[y_1] \ss \ol{S}$ is finite and integral. Then $y_0 \in \NZD(\ol{R})$, hence $y_0$ is a multiplicity parameter of degree $1$ for $\ol{R}$. Therefore \[e_0(\ol{S}) = e_0(\ol{R}/y_0\ol{R}) = e_0(\ol{R}).\] 
Now, let $\ol{h} \in \ol{R}$ be a homogeneous generator of $\ol{\a}\sat$. According to Lemma \ref{asatGen} we have $\ol{h} = h_0y_0^{k_0+1}+\ol{g} \in R_{k_0+1}$ with $h_0 \in K\bs\{0\}$ and $\ol{g} \in \ol\p$, and as $y_0 \in \NZD(\ol{R})$ it follows $\ol{h} \in \NZD(\ol{R})$, so $\ol{h}$ is a multiplicity parameter of degree $k_0+1$ for $\ol{R}$. As a result, we finally get \[\setlength\arraycolsep{0.2em}\begin{array}{lclclcl} e_0(R/(\a+\q R)) &=& e_0(\ol{R}/\ol{\a}) &=& e_0(\ol{R}/\ol{\a}\sat)\\ &=& e_0(\ol{R}/\ol{h}\ol{R}) &=& (k_0+1)\cdot e_0(\ol{R}) &=& (k_0+1)\cdot e_0(\ol{S}).\end{array}\]
\end{proof}

As a corollary to Theorem \ref{Theorem} we get the main result about PEIs (see \cite[Theorem 3.5]{C-S}, \cite[Proposition 6.2]{G}):

\begin{cor}\label{Green} Assume $\a \nss \p$, and let $\q \in \mProj(S)$. Then \[\K_k^\p(\a) \ss \q \Llra e_0(R/(\a + \q R)) > k.\]
\end{cor}

\begin{proof} We use the same notations as above. If $\K_0^\p(\a) = \a\cap S \nss \q$, then $\sqrt{\a\cap S + \q} = S_+$ and $e_0(\ol{S}) = 0$. So, assume $\a \cap S \ss \q$. Then $\ol{S} = S/\q = K[y_1]$, so $\dim(\ol{S}) = 1 = e_0(\ol{S})$, and $\ol{\a} \ss \ol{R} = K[y_0,y_1]$ is an ideal of height $1$. Therefore $\ol{\a}\sat$ is a principal ideal, and we get our claim by Theorem \ref{Theorem}.	 
\end{proof}

\begin{r&n}\label{Notat} (A) If $Z \ss \PP^n$ is a closed subscheme, by {\it the homogeneous ideal of $Z$} we mean the unique saturated ideal $I_Z \ss R$ such that $Z = \Proj(R/I_Z)$ (see \cite[II, Exercise 5.10]{H}). 
\\
(B) If $X,Y \ss \PP^n$ are closed subschemes given by the homogeneous ideals $I_X, I_Y \ss R$ such that $\dim(X\cap Y) = 0$, we denote by \[l(X\cap Y) = e_0(R/(I_X+I_Y))\] the {\it length of the intersection of $X$ and $Y$}.
\\
(C) Let $\Lambda, \Delta \ss \PP^n$ be linear subspaces with homogeneous ideals $I_{\Lambda}$ and $I_\Delta \ss R$, respectively; these ideals are generated by linear forms. The {\it linear span $\langle \Lambda, \Delta\rangle$ of $\Lambda$ and $\Delta$} is defined as the linear subspace defined by the common linear forms of $I_\Lambda$ and $I_\Delta$, that is \[\langle\Lambda,\Delta\rangle = \Proj(R/((I_\Lambda)_1\cap (I_{\Delta})_1)R).\]
\\
(D) Let $Z \ss \PP^n$ be a closed subscheme, let $p \in \PP^n\bs Z$, and let $k \in \Nat$. A $k$-secant line to $Z$ is a line $\LL \ss \PP^n$ such that $l(Z \cap \LL) \geq k$. We define the {\it $k$-secant cone $\Sec_p^k(Z)$ of $Z$ with vertex $p$} as the closed subset of $\PP^n$ \[\Sec^k_p(Z) \defin \{p\}\cup \bigcup\{\LL \mid \LL \mbox{ is a } k \mbox{-secant line to } Z \mbox{ with } p \in \LL\}\] furnished with its structure of reduced closed subscheme of $\PP^n$. Next, we define the {\it $k$-secant loci of $Z$ with respect to $p$} as the closed subscheme of $\PP^n$ \[\Sigma^k_p(Z) \defin Z \cap \Sec_p^k(Z).\] Some authors also use the the term ``entry locus'' instead of ``secant locus''. Observe that $\Sec_p^k(Z) = \Jn(p, \Sigma_p^k(Z))$ is the (embedded) join of $p$ and $\Sigma_p^k(Z)$; this is the reason we demand $\Sec_p^k(Z)$ to be reduced (see \cite{FCV}). Example \ref{sat} shows the importance of defining the secant cone to be reduced.
\\
(E) Let $p \in \PP^n$. We denote the linear projection with centre in $p$ by \[\pi_p: \PP^n_K\bs\{p\} \lra \PP^{n-1} = \PP^{n-1}_K(p);\] it is given by $S = K[p_1]\hookrightarrow R$. Let $\pi_p\res_Z$ denote the restriction of $\pi_p$ to a closed subscheme $Z \ss \PP^n$ with homogeneous ideal $I_Z$. The {\it fibre of $\pi_p\res_Z$ over a closed point $q \in \pi_p(Z)$} is \[(\pi_p\res_Z)^{-1}(q) = \langle q,p\rangle\cap Z.\] As $q \in \mProj(K[p_1])$ the fibre $(\pi_p\res_Z)^{-1}(q)$ is given by the ideal $(qR+I_Z)\sat\ss R$.
\end{r&n}

\begin{rem} (A) Keep the notations of above. Now, Corollary \ref{Green} just says that if $p \notin Z$, then the $k$-th PEI $\K^p_k(I_Z)$ defines the set of closed points $q \in \pi_p(Z)$ whose fibres $(\pi_p\res_Z)^{-1}(q)$ are of length $>k$. Note that $\K^p_k(I_Z)$ does not need to be saturated (see Example \ref{sat}).\\
(B) In this paper, we only consider outer projections, that is projections from a point $p$ not contained in $Z$, because we are mainly interested in secant cones. But the study of PEIs can be useful for the study of projections from $p$ if $p \in Z$; see \cite{H-K} for an application of PEIs to inner projections.
\end{rem}

\begin{thm}\label{SecPEI} Let $Z \ss \PP^n$ be a closed subscheme with homogeneous ideal $I_Z \ss R$, and let $p \in \PP^n$ such that $p \notin Z$. Let $k \in \Nat$. Then \[\Sec^k_p(Z) = \Proj\left(R/\sqrt{\K^p_{k-1}(I_Z)R}\right)\] and \[\Sigma^k_p(Z) = \Proj\left(R/\left(I_Z + \sqrt{\K^p_{k-1}(I_Z)R}\right)\right).\]
\end{thm}

\begin{proof} Let $q \in \PP^{n-1} = \mProj(S)$; the homogeneous ideal of the projective line $\langle q,p\rangle \ss \PP^n$ is $qR \in \Proj(R)$. According to corollary \ref{Green} it therefore holds \[\setlength\arraycolsep{0.2em}\begin{array}{rclcl} V_{\PP^{n-1}}(\K^p_{k-1}(I_Z)) &=& \{q \in \pi_p(Z) \mid e_0(R/(I_Z + qR)) \geq k\}\\ &=& \{q \in \pi_p(Z) \mid l(Z \cap \langle p, q\rangle) \geq k\} &=& \pi_p(\Sigma^k_p(Z)).\end{array}\] 
But as the closure of $\pi_p^{-1}(\pi_p(\Sigma^k_p(Z))) \ss \PP^n$ is just $\Sec^k_p(Z)$, we get the first equation. The second equations follows by definition.
\end{proof}

\section{Multiple Projections}

Let $\widetilde{S} \ss R$ be a homogeneous graded $K$-subalgebra, that is there exist an integer $t \in \{0, \ldots,n\}$ and linearly independent elements $y_0, \ldots, y_t \in R_1$ such that $R = \widetilde{S}[y_0,\ldots,y_t]$. For a graded ideal $\a \ss R$, let $\tilde{\a} \defin \a\cap\widetilde{S}$. Note that there is a natural inclusion map $\widetilde{S}/\widetilde{\a} \inkl R/\a$; we therefore consider $\widetilde{S}/\widetilde{\a}$ as a graded $K$-subalgebra of $R/\a$.

\begin{thm}\label{MultProjAlg} 
 Let $\a \ss R$ be a graded ideal such that $(R/\a)_m = (\widetilde{S}/\widetilde\a)_m$ for all $m\>>0$. Let $\widetilde{\p} \in \mProj(\widetilde{S})$ such that $\widetilde\a \nss \widetilde{\p}$, let $S \defin K[\widetilde{\p}_1] \ss \widetilde{S}$, and let $\q \ss S$ be a graded ideal such that $\dim(S/(\widetilde{\a}\cap S + \q))=1$. Let $l_0 \defin \max\{l \in \Nat\cup\{-1\} \mid \K^{\widetilde{\p}}_l(\widetilde{\a})\ss\widetilde{\a}\cap S + \q\}$. Then \[e_0(R/(\a+\q R)) = (l_0+1)\cdot e_0(S/(\a\cap S+\q)).\]
\end{thm}

\begin{proof}
 Let $m \>> 0$. As $(R/\a)_m = (\widetilde{S}/\widetilde{\a})_m$ it holds $((y_0, \ldots, y_t)R)_{m} \ss \a_m$ and in particular $((y_0, \ldots,y_t)\q)_m \ss \a_m$. Moreover $\q R = \q\widetilde{S}+(y_0,\ldots,y_t)\q$. It follows
 \[\setlength\arraycolsep{0.2em}\begin{array}{rclcl} (R/(\a+\q R))_m &=& (R/(\a+\q\widetilde{S}))_m\\ &=&
 ((\widetilde{S}/\widetilde{\a})/(\widetilde{\a}+\q\widetilde{S}/\widetilde{\a}))_m &=& (\widetilde{S}/(\widetilde{\a}+\q\widetilde{S}))_m\end{array}\]
 Thus, by \ref{Theorem} \[\setlength\arraycolsep{0.2em}\begin{array}{lclcl} e_0(R/(\a+\q R)) &=& e_0(\widetilde{S}/(\widetilde{\a}+\q\widetilde{S}))\\ &=& (l_0+1)\cdot e_0(S/(\widetilde{\a}\cap S+\q)) = (l_0+1)\cdot e_0(S/(\a\cap S+\q)).\end{array}\]
\end{proof}

Now, let $d \in \Nat$, let $\Lambda = \PP^d \ss \PP^n$ and $S \defin K[\widetilde{\p}_1] \ss \widetilde{S}$ a linear subspace of dimension $d$ with homogeneous ideal $\I_\Lambda\ss R$, and let \[\pi \defin \pi_\Lambda: \PP^n\bs\Lambda \lra \PP^{n-d-1}\] be the projection with centre in $\Lambda$; this projection is given by $S = S_\Lambda \defin K[(I_\Lambda)_1] \inkl R$. Then, we can choose a decomposition \[\pi = \pi_d\circ\pi_{d-1}\circ\cdots\circ\pi_0,\] where $\pi_i:\PP_K^{n-i}(p_{i-1})\bs\{p_i\}\lra\PP^{n-i-1}_K(p_i)$ are simple projections for $i \in \{0,\ldots,d\}$. If we denote the homogeneous rings of $\PP^{n-1}(p_0),\ldots,\PP^{n-d}(p_{d-1})$ by $\widetilde{S}_0, \ldots,\widetilde{S}_{d-1}$, this decomposition is given by \[S \inkl \widetilde{S}_{d-1} \inkl \cdots \widetilde{S}_0 \inkl R.\]

\begin{cor}\label{SemiIsoProj}
 Let $Z \ss \PP^n$ be a closed subscheme with homogeneous ideal $I_Z$, and assume that there is a decomposition $\pi = \pi_d\circ\cdots\circ\pi_0$ such that $(\pi_{d-1}\circ\cdots\circ\pi_0)\res_Z: Z \to \pi_d^{-1}(\pi(Z))$ is an isomorphism. Let $q \in \PP^{n-d-1}$ be a closed point. Then, for all $k \in \Nat$ \[l(Z\cap\langle q,\Lambda\rangle) > k \Llra \K^{I_\Lambda\cap \widetilde{S}_{d-1}}_k(I_Z\cap \widetilde{S}_{d-1}) \ss \q.\]
\end{cor}

\begin{proof}
 $(\pi_{d-1}\circ\cdots\pi_0)\res_Z$ being an isoprojection is equivalently to $(R/I_Z)_m = (\widetilde{S}_{d-1}/(I_Z\cap \widetilde{S}_{d-1}))_m$ for all $m \>> 0$. So, we get our claim by Proposition \ref{MultProjAlg}.
\end{proof}

\begin{notat} 
For the remainder of this section, let $\L \in \Proj(R)$ be a linearly generated ideal of height $n-1$, that is $\LL \defin \Proj(R/\L) = \PP^1 \ss \PP^n$ is a projective line. Let $S = S_\L \defin K[\L_1] \ss R$; the twofold projection $\pi_{\LL}: \PP^n\bs\LL \Epi \PP^{n-2}$ is given by $S \inkl R$. Let $p,p' \in \LL, p\neq p'$, and let $\widetilde{S} \defin K[p_1],\widetilde{S}' \defin K[p'_1] \ss R$. Consider the projections $\pi: \PP^n\bs\{p\} \lra \PP^{n-1}$ and $\pi': \PP^n\bs\{p'\} \lra \PP^{n-1}$ given by $\widetilde{S} \inkl R$ and $\widetilde{S}' \inkl R$, respectively, as well as $\tilde{\pi}: \PP^{n-1}\bs\{\tilde{p}\} \lra \PP^{n-2}$ and $\tilde{\pi}': \PP^{n-1}\bs\{\tilde{p}'\} \lra \PP^{n-2}$ given by $S \inkl \widetilde{S}$ and $S \inkl \widetilde{S}'$, respectively, where $\tilde{p} \defin \pi'(p) = \pi'(\LL\bs\{p'\})$ and $\tilde{p}' \defin \pi(p') = \pi(\LL\bs\{p\})$. Then \[\pi_{\LL} = \tilde{\pi}'\circ\pi = \tilde{\pi}\circ\pi': \PP^n\bs\LL \lra \PP^{n-2}\] are two decompositions of $\pi$.
\end{notat}

\begin{definition}\label{ClDecomp}
 Let $\a \ss R$ be a graded ideal such that $R_+ \ss \sqrt{\a+\L}$. We call $\p,\p' \in \mProj(R)\cap\Var(\L)$ a {\it clever decomposition of $\L$ with respect to $\a$} if \[\a\sat = \left((\a\cap K[\p_1])R+(\a\cap K[\p'_1])R\right)\sat.\]
\end{definition}

\begin{rem}
 Let $Z \ss \PP^n$ be a closed subscheme with homogeneous ideal $I_Z \ss R$ such that $Z\cap \LL = \emptyset$. Geometrically, Definition \ref{ClDecomp} means that $p,p' \in \LL$ are a {\it clever decomposition of $\LL$ with respect to $Z$} if \[Z = \Jn(\pi(Z),p) \cap \Jn(\pi'(Z),p') = \Jn(Z,p) \cap \Jn(Z,p').\]
\end{rem}

\begin{thm}\label{ClMultProjAlg}
 Let $\a \ss R$ be a graded ideal such that $R_+ \ss \sqrt{\a+\L}$, and let $\p,\p' \in \mProj(R)\cap\Var(\L)$ be a clever decomposition of $\L$ with respect to $\a$. Let $\q \in \mProj(S)$, and let $k_0 \defin \max\{k \in \Nat\cup\{-1\} \mid \K_k^{\L\cap \widetilde{S}}(\a\cap \widetilde{S}) \ss \q\}$, $k'_0 \defin \max\{k \in \Nat\cup\{-1\} \mid \K_k^{\L\cap \widetilde{S}'}(\a\cap \widetilde{S}') \ss \q\}$. Then \[e_0(R/(\a+\q R)) = (k_0+1)\cdot(k'_0+1).\]
\end{thm}

\begin{proof}
 We write $\ol{R} \defin R/\q R$. Let $y_i0 \in R_1$ such that $\p = \L+y_0R$, and let $y_2 \in S_1\ss R_1$ such that $\q +y_2R = \L$, that is $\widetilde{S}/\q \widetilde{S} = K[y_0,y_2]$ and $\ol{R} = K[y_0,y_1,y_2]$. By Lemma \ref{asatGen}, there is a homogeneous element 
 \[\ol{h} = h_0y_0^{k_0+1}+h_1y_0^{k_0}y_2+\cdots+h_{k_0+1}y_2^{k_0+1} \in (K[y_0,y_2])_{k_0+1} \ss \ol{R}_{k_0+1},\] where $h_0, \ldots, h_{k_0+1} \in K$ such that $h_0 \neq 0$ and \[\a\cap \widetilde{S} +\q \widetilde{S}/\q \widetilde{S} = \ol{h}\ol{S}.\] Analogously, we can choose $y'_0 \in R_1\bs\p'_1$ such that $\a\cap \widetilde{S} +\q \widetilde{S}/\q \widetilde{S}$ is generated by \[\ol{h}' = h'_0y_0'^{k_0+1}+h'_1{y'}_0^{k'_0}y_2+\cdots+h'_{k'_0+1}y_2^{k'_0+1} \in (K[y'_0,y_2])_{k'_0+1} \ss \ol{R}_{k'_0+1},\] where $h,_0, \ldots, h,_{k,_0+1} \in K$ such that $h'_0 \neq 0$
 As $h_0,h'_0 \neq 0$ it follows that $\ol{h},\ol{h}' \ss \ol{R} = K[y_0,y'_0,y_2]$ is an $\ol{R}$-regular sequence. So, $\ol{h},\ol{h}'$ is a system of multiplicity parameters of degree $k_0+1,k'_0+1$ for $\ol{R}$. 
 Further, for all $d \>> 0$ \[\setlength\arraycolsep{0.2em}\begin{array}{rclclcl} (\a+\q R/\q R)_d &=& \a\sat_d + (\q R)_d/(\q R)_d\\ &=& \left((\a\cap \widetilde{S})R + (\a\cap \widetilde{S}')R\right)_d + (\q R)_d/(\q R)_d\\ &=& \left((\a\cap \widetilde{S} +\q \widetilde{S}/\q \widetilde{S})\ol{R}\right)_d + \left((\a\cap \widetilde{S}' +\q \widetilde{S}'/\q \widetilde{S}')\ol{R}\right)_d\\ &=& \left((\ol{h},\ol{h}')\ol{R}\right)_d.\end{array}\] 
 Therefore, we get \[e_0(R/(\a+\q R)) = e_0(\ol{R}/(\ol{h},\ol{h}')\ol{R}) = (k_0+1)\cdot(k'_0+1).\]
\end{proof}

\begin{cor}\label{ClMultProj}
 Let $Z \ss \PP^n$ be a closed subscheme such that $Z\cap \LL = \emptyset$, and let $p,p' \in \LL$ be a clever decompositon of $\LL$ with respect to $Z$. Then, for all closed points $q \in \PP^{n-2}$ \[l(Z\cap\langle q,\LL\rangle) = l(\pi(Z)\cap\langle q,\pi(p')\rangle)\cdot l(\pi'(Z)\cap\langle q,\pi'(p)\rangle).\]
\end{cor}

\begin{proof}
 Clear by Proposition \ref{ClMultProjAlg} and Theorem \ref{Theorem}.
\end{proof}

\begin{rem}
 The next obvious questions here would be whether there is a clever decomposition for any $\L$ and any $\a$, and if not what conditions on $\L$ and $\a$ imply the existence of a clever decomposition. We are not going to answer these questions here; for now, we are not interested in clever decompositions themselves but in their usefulness for studying examples (see Example \ref{ClDecEx}).
\end{rem}

\section{Computational Aspects}

We keep the previous notations.

\begin{thm}\label{Groebner} Let $\a \ss R$ be a graded ideal, and assume $\p \defin (x_1, \ldots, x_n) \in \mProj(R)$. Let $\sigma$ be an elimination ordering on $R$, and let $G$ be a Gr{\"o}bner basis of $\a$ with respect to $\sigma$. For all $k \in \Nat$, the set \[G_k\defin \{\LC_{x_0}(g) \mid g \in G \wedge \deg_{x_0}(g) \leq k\}\] is a Gr{\"o}bner basis of $\K_k^{\p}(\a)$ with respect to the term ordering on $S = K[\p_1]$ induced by $\sigma$.
\end{thm}

\begin{proof} \cite[Proposition 3.4]{C-S} and Corollary \ref{JustOneCoord}.
\end{proof}

\begin{rem}\label{k0}
 We keep the notations of Proposition \ref{Groebner}, but we assume further that $\a \nss \p$, that is $R_+ \ss \sqrt{\a+\p}$. So, there is an integer $t \in \Nat$ such that $x_0^t \in \a+\p$, hence $x_0^t$ is contained in the initial ideal $\In_\sigma(\a+\p)$. But as $\p = (x_1, \ldots, x_n)$, this means $x_0^t \in \In_\sigma(\a)$. Therefore, there must be an element $g_0 \in G$ such that $\LT_{x_0}(g_0) = x_0^s$ for some $s \leq t$ and $\K^{\p}_s(\a) = K[\p_1]$. 
\end{rem}

\begin{algor}\label{SecComp} (A) Using Propostion \ref{Groebner}, we obtain the following method for computing the ideals of secant cones and secant loci:
\medskip\\
Let $Z \ss \PP^n$ be a closed subscheme with homogeneous ideal $I_Z \ss R$, and let $p \in \PP^n$ such that $p \notin Z$. 
\\
1. Choose a linear coordinate transformation $\psi: R \stackrel{\cong}{\to} R$ such that $\psi(p) = (x_1, \ldots, x_n)$.
\\
2. Compute a Gr{\"o}bner basis $G$ of $\psi(I_Z)$ with respect to an elimination ordering.
\\
3. Choose $k_0 \in \Nat$ such that $\K^{\psi(p)}_k(\psi(I_Z)) = K[p_1]$ for all $k \geq k_0$. An integer $k_0$ with this property exists by Remark \ref{k0}.
\\
4. Compute the partial elimination ideals $\K^{\psi(p)}_0(\psi(I_Z)), \ldots, \K^{\psi(p)}_{k_0-1}(\psi(I_Z))$. This can easily be done using Proposition \ref{Groebner}.
\\
5. Set $\K^p_k(I_Z) \defin \psi^{-1}(\K^{\psi(p)}_k(\psi(I_Z)))$ for $0 \leq k \leq k_0-1$. Lemma \ref{CoordTransf} guarantees that we indeed get the partial elimination ideals of $I_Z$ with respect to $p$. 
\\
6. Compute $\sqrt{\K^p_k(I_Z)R}$ for $0 \leq k \leq k_0-1$.
\\
By Proposition \ref{SecPEI}, for any $k \in \{0, \ldots, k_0-1\}$, the $(k+1)$-secant cone of $Z$ with respect to $p$ is defined as a scheme by the homogeneous ideal $\sqrt{\K_k^p(I_Z)R}$, while the $(k+1)$-secant loci of $Z$ with respect to $p$ is defined by the homogeneous ideal $\sqrt{\K_k^p(I_Z)R}+I_Z$. As $\K_k^p(I_Z)R = R$ for all $k \geq k_0$, the higher secant loci $\Sigma^k_p(Z)$ are empty.
\medskip\\
(B) The above method contains some choices. We can replace these choices with explicit terms and get the following algorithm:
\medskip\\
{\bf Input:} The homogeneous ideal $I_Z \ss R$ of a closed subscheme $Z \ss \PP^n$, and a minimal system of generators $y_1, \ldots, y_n \in R_1$ of the closed point $p \in \Proj(R)$. Consider $R$ to be furnished with either the lexicographical term order or the reversed lexicographical term order.
\\
1.1. Compute $l \defin \min\{i \in \{0, \ldots,n\} \mid x_i \notin p\}$. 
\\
1.2. Define the coordinate transformation $\psi: R \stackrel{\cong}{\lra} R$ to be the inverse of $x_0 \mapsto x_l$, $x_1 \mapsto y_1, \ldots, x_n \mapsto y_n$. Then indeed $\psi(p) = (x_1, \ldots, x_n)$. Calculate $\psi(I_Z)$. 
\\
2. Compute a Gr{\"o}bner basis $G$ of $\psi(I_Z)$, for example using the Buchberger algorithm.
\\
3. Set $k_0 \defin \max\{\deg_{x_0}(g) \mid g \in G\}$. Then $\K^{\psi(p)}_k(\psi(I_Z)) = K[p_1]$ for all $k \geq k_0$ according to Remark \ref{k0}.
\\
4. For all $k \in \{0, \ldots, k_0-1\}$, set $G_k \defin \{\LT_{x_0}(g) \mid g \in G \wedge \deg_{x_0}(g) \leq k\}$.
\\
5. Set $\K^p_k(I_Z)R \defin \psi^{-1}(G_k)R$ for $0 \leq k \leq k_0-1$. 
\\
6. Compute $\sqrt{\K^p_k(I_Z)R}$ for $0 \leq k \leq k_0-1$, for example using the algorithm of Krick and Logar.
\\
{\bf Output:} Ideals $\sqrt{\K_0^p(I_Z)R}, \ldots, \sqrt{\K_{k_0-1}^p(I_Z)R}$ (via a finite set of generators) of $\Sec^p_1(Z), \ldots, \Sec^p_{k_0}(Z)$.
\end{algor}

An implementation of this algorithm for {\sc Singular} can be obtained on request from the author (even if it is still rather unpolished).

\section{Examples}

We use the notations of the previous sections.

\begin{ex} (A) Let $R \defin K[x_0,\ldots,x_4]$, let $\p \defin (x_1,\ldots,x_4)$, let \[\a \defin (x_0^4+x_1^2x_2^2,x_0^2x_1-x_3^3,x_2^2-x_3^2,x_0x_2+x_4^2) \ss R,\] and let $\q \defin (x_3,x_4) \ss S = K[x_1,\ldots,x_4]$. Then $\a\cap S +\q = (x_2^2,x_3,x_4)$ and $x_2 \in \K^\p_1(\a)$, so $k_0 = 0$ and $e_0(\ol{S}) = e_0(K[x_1,\ldots,x_4]/(x_2^2,x_3,x_4)) = 2$. Furthermore, \[\a + \q R = (x_0^4,x_0^2x_1,x_0x_2,x_2^2,x_3,x_4)\] and therefore \[e_0(R/(\a+\q R)) = 3 > (k_0+1)\cdot e_0(S/(\a\cap S + \q)).\]
(B) Keep $R$, $\p$, $S$ and $\q$ of part (A), and let \[\a \defin (x_0^4+x_0x_1^3,x_0^3x_1+x_1^4+x_3^4, x_0^2x_2+x_4^3,x_2^2) \ss R.\] Then, \[\K_1^\p(\a) = \a\cap S +(x_3^4,x_1x_4^3) \ss \a\cap S + \q = (x^2_2,x_3,x_4),\] but $x_2 \in \K^\p_2(\a)$, so $k_0 = 1$ and $e_0(\ol{S}) = 2$. On the other hand, \[\a + \q R = (x_0^4+x_0x_1^3,x_0^3x_1+x_1^4, x_0^2x_2,x_0x_1^3x_2,x_1^4x_2,x_2^2,x_3,x_4)\] and hence \[e_0(R/(\a+\q R)) = 3 < (k_0+1)\cdot e_0(S/(\a\cap S +q)).\]
(C) Again, keep $R$, $\p$ and $S$ as in part (A), but now let \[\a \defin (x_0^5+x_0^2x_1^3,x_0^4x_1+x_0x_1^4+x_3^5,x_0^3x_1^2+x_1^5+x_4^5,x_0^3x_2,x_2^3),\] and let $\q \defin (x_1x_2,x_3,x_4) \ss S$. Then, \[\K_2^\p(\a) = \a\cap S + (x_1^3x_2,x_4^5,x_3^5) \ss \a\cap S + \q = (x_1x_2,x_2^3,x_2,x_4),\] but $x_2 \in \K_3^\p(\a)$, so $k_0 = 2$ and $e_0(\ol{S}) = 1$. On the other hand, \[\ol{\a}\sat = (x_0^3+x_1^3,\ol{x_2}) \ss \ol{R} = K[x_0,x_1,\ol{x_2}]\] and therefore \[e_0(R/(\a+\q R)) = 3 = (k_0+1)\cdot e_0(S/(\a\cap S+\q)).\]
(D) Now, let $R = K[x_0,\ldots,x_n]$, $\p \in \mProj(R)$, let $\a \ss R$ be a graded ideal such that $\a \nss \p$, and let $\q \ss S\defin K[\p_1]$ be a graded ideal such that $\dim(\ol{S}) = 1$. Then, $\ol{\a}_0 \defin (\ol{\a}\cap K[y_0,y_1])\sat$ is a principal ideal; the same argument as used in the proof of Lemma \ref{asatGen} shows that a homogeneous generator $\ol{h}$ of $\ol{\a}_0$ is of degree $\geq k_0+1$. As $\ol{\a}_0\ol{R} \ss \ol{\a}\sat$, it follows \[e_0(R/(\a+\q R)) \leq e_0(\ol{R}/\ol{\a}_0\ol{R}) = e_0(\ol{R}/\ol{h}\ol{R}) \geq (k_0+1)e_0(S/(\a\cap S + \q)).\] Now, (A) and (B) above prove that both inequalities between $e_0(R/(\a+\q R))$ and $(k_0+1)\cdot e_0(S/(\a\cap S + \q))$ can occur if $\ol{\a}\sat$ is not a principal ideal. But the condition that $\ol{\a}$ is a principal ideal is not necessary for $e_0(R/(\a+\q R)) = (k_0+1)\cdot e_0(S/(\a\cap S + \q))$, as (C) shows. Indeed, we conjecture that this equality always holds if $\ol{S}_+$ can be generated by two elements and the degree of a generator of $\ol{\a}_0$ is $k_0+1$.
\end{ex}

\begin{ex} Let $X \ss \PP^n$ be a smooth rational normal scroll and of codimension at least $2$. Let $p \in \PP^n\bs X$ be a closed point. Then, according to \cite[Theorem 3.2]{B-P}, either 
\begin{enumerate}
 \item[(a)] $\Sec_p^2(X) = \PP^1$ and $\Sigma_p^2(X)\ss \PP^1$ is either a double point or the union of two simple points.
 \item[(b)] $\Sec_p^2(X) = \PP^2$ and $\Sigma_p^2(X)\ss \PP^2$ is either a smooth conic or the union of two lines $L, L' \ss X$. 
 \item[(c)] $\Sec_p^2(X) = \PP^3$ and $\Sigma_p^2(X)\ss \PP^3$ is a smooth quadric surface.
 \item[(d)] $\Sec_p^2(X) = \emptyset$, i.e. $p \notin \Sec(X)$.
\end{enumerate}
Now, let us consider the scroll $X = S(1,1,2,3) \ss \PP^{10}$. Using Algorithm \ref{SecComp}, it is easy to find 6 points of $\PP^{10}\bs S(1,1,2,3)$ such that every one of the possible six secant loci occurs (see Table 1; there $\mathfrak{A}_p$ of a closed point $p = (p_0: \cdots: p_{10}) \in \PP^{10}\bs X$ is the set of indices $i \in \{0,\ldots,10\}$ such that $p_i = 1$ for $i \in \mathfrak{A}_p$ and $p_i = 0$ else).
\begin{table}\centering\begin{tabular}{|c|c|c|c|c|}\hline
 $\mathfrak{A}_p$ & $\K^p_1(I_X)$ & $\Sec_p^2(X)$ & $\begin{array}{c}\K^p_1(I_X)R\\+I_X/K^p_1(I_X)R\end{array}$ & $\Sigma_p^2(X)$ \\ \hline
 $\{7,10\}$ & $\left(\begin{array}{c}x_0,\ldots,x_6,\\x_8,x_9\end{array}\right)$ & $\PP^1$ & $\begin{array}{c}(x_7x_{10})\\ \ss K[x_7,x_{10}]\end{array}$ & $\begin{array}{c}\mbox{Two} \\ \mbox{points}\end{array}$\\ \hline
 $\{9\}$ & $(x_0,\ldots,x_8)$ & $\PP^1$ & $\begin{array}{c}(x_9^2)\\ \ss K[x_9,x_{10}]\end{array}$ & $\begin{array}{c}\mbox{Double} \\ \mbox{point}\end{array}$\\ \hline
 $\{0,10\}$ & $(x_2,\ldots,x_9)$ & $\PP^2$ & $\begin{array}{c}(x_0x_{10})\\ \ss K[x_0,x_1,x_{10}]\end{array}$ & $\begin{array}{c}\mbox{Two} \\ \mbox{lines}\end{array}$\\ \hline
 $\{5\}$ & $\left(\begin{array}{c}x_0,\ldots,x_3,\\x_7,\ldots,x_{10}\end{array}\right)$ & $\PP^2$ & $\begin{array}{c}(x_4x_6-x_5^2)\\ \ss K[x_4,x_5,x_6]\end{array}$ & $\begin{array}{c}\mbox{Smooth} \\ \mbox{conic}\end{array}$\\ \hline
 $\{0,3\}$ & $\left(\begin{array}{c}x_4,x_5,x_6,\\x_8,x_9,x_{10}\end{array}\right)$ & $\PP^3$ & $\begin{array}{c}(x_0x_3-x_1x_2)\\ \ss K[x_0,\ldots,x_3]\end{array}$ & $\begin{array}{c}\mbox{Quadric} \\ \mbox{surface}\end{array}$\\ \hline
 $\{4,9\}$ & $\left(\begin{array}{c}x_0,\ldots,x_3,\\x_5,\ldots, x_8,\\x_{10},x_4-x_9\end{array}\right)$ & $p$ & $(x_9^2)\ss K[x_9]$ & $\emptyset$\\ \hline
\end{tabular} \caption{Examples of all possible secant loci of $X = S(1,1,2,3)$}\end{table}
\end{ex}

\begin{ex}\label{sat}
Let $Z \ss \PP^3$ be the subscheme defined by the homogeneous ideal \[I_Z \defin (x_0^4,x_0^3x_1,x_0^2x_1^2+x_0x_3^3,x_0x_1x_2^2+x_1^4,x_1x_3^3) \ss R \defin K[x_0,x_1,x_2,x_3].\] As $\sqrt{I_Z} = (x_0,x_1)$, the underlying set of $Z$ is just a line, and $I_Z$ is saturated. Computing the PEIs of $I_Z$ with respect to the point $p = (1:0:0:0) \in \PP^3$ we get \[\setlength\arraycolsep{0.2em}\begin{array}{lcl} \K^p_0(I_Z) = (x_1x_3^3,x_1^9) &\ss& S \defin K[x_1,x_2,x_3],\\ \K_1^p(I_Z) = (x_1x_3^3,x_1^5-x_2^2x_3^3,x_1x_2^2,x_3^6),\\ \K_2^p(I_Z) = (x_1^2,x_1x_2^2,x_3^3),\\ \K_3^p(I_Z) = (x_1,x_3^3),\\ \K_4^p(I_Z) = S.\end{array}\]
So, we can compute $\K^p_1(I_Z)\sat = \K^p_2(I_Z)\sat = \K_3^p(I_Z)$, that is the first and second PEI of $I_Z$ with respect to $p$ are not saturated; their saturation indices are $8$ and $4$, respectively. Furthermore, we see that $\Sigma_p^2(Z) = \Sigma_p^3(Z) = \Sigma_p^4(Z)$ are equal as sets and consist just of the point $q = (0:0:1:0)$. The line $\langle q,p\rangle$ is a $4$-secant to $Z$. Therefore, $Z$ is not smooth.\\
On the other hand, the extension ideals $\K^p_1(I_Z)R$, $\K^p_2(I_Z)R$ and $\K^p_3(I_Z)R$ are saturated, meaning that the schemes $\Proj(R/\K^p_1(I_Z)R),$ $\Proj(R/\K^p_2(I_Z)R)$ and $\Proj(R/\K^p_3(I_Z)R)$ are different; they are non-reduced and therefore not equal to the $k$-secant cones of $Z$ for $k \in \{2,3,4\}$. Another consultation of {\sc Singular} tells us that \[(I_Z+\K^p_1(I_Z)R)\sat = (I_Z+\K^p_2(I_Z)R)\sat = (I_Z+\K^p_3(I_Z)R)\sat = (x_0^4,x_1,x_3^3),\] while \[\begin{array}{rclcl} \left(I_Z+\sqrt{\K^p_1(I_Z)R}\right)\sat &=& \left(I_Z+\sqrt{\K^p_2(I_Z)R}\right)\sat\\ &=& \left(I_Z+\sqrt{\K^p_3(I_Z)R}\right)\sat &=& (x_0^4,x_1,x_3).\end{array}\]
Hence, we indeed have to demand that $\Sec_p^k(Z)$ is reduced; if we omitted this condition, we would get $l(Z \cap \Sec_p^4(Z)) = 12$, where $\Sec_p^4(Z) = \langle p,q\rangle$ is just a line. But $\langle p,q\rangle$ certainly is no $12$-secant line to $Z$. 
\end{ex}

\begin{ex}[Example 7.4(E) in \cite{B-S}]\label{7.4(E)}
Let $n = 10$. Consider the rational normal scroll $W \defin S(1,8) \ss \PP^{10}_K$ with homogeneous ideal $I_W\ss R$, and let $\LL\ss\PP^{10}_K$ be the line given by $\L = (x_0,x_1,x_2,x_5,\ldots,x_{10})\ss R$. Let $\pi_{\LL}:\PP^{10}_K\bs\LL \lra \PP^8_K$ be the double projection given by $S \defin K[x_0,x_1,x_2,x_5,\ldots,x_{10}] \inkl R$, and let $Y \defin \pi_{\LL}(W) \ss \PP^8_K$. The homogeneous ideal $J \ss S$ of $Y$ is given by $18$ quadrics and one quartic $Q = x_1^3x_2-x_0^3x_5$. 
Now, let us consider the secant loci of $W$ with respect to the points of $\LL$. According to \cite{CJ}, the secant variety of $W$ is given by the ideal $M$ generated by the $3\times3$-minors of the matrix 
\[\left(\begin{array}{ccccccc} x_2 & x_3 & x_4 & x_5 & x_6 & x_7 & x_8\\ x_3 & x_4 & x_5 & x_6 & x_7 & x_8 & x_9\\ x_4 & x_5 & x_6 & x_7 & x_8 & x_9 & x_{10}\end{array}\right),\]
so that $M+\L = (x_0,x_1,x_2,x_4^3,x_5,\ldots,x_{10}) \ss R$ and therefore $\Sec(W)\cap\LL$ contains just one point $p = (0:0:0:1:0:\cdots:0)$ with homogeneous ideal $(x_0,x_1,x_2,x_4,\ldots,x_{10}) \in \mProj(R)$. The partial elimination ideals of $I_W$ with respect to $p$ are 
\[\setlength\arraycolsep{0.2em}\begin{array}{l}\K^{p}_0(I_W) = I_W\cap K[(p_0)_1],\\ \K^{p}_1(I_W) = (x_0,x_1,x_4,\ldots,x_{10}),\\ \K^{p}_2(I_W) = K[(p_0)_1],\end{array}\]
so that \[\Sec^2_{p_0}(W) = V_{\PP^{10}_K}(\K_1^{p_0}(I_W)R) = \PP^1 \ss \PP^{10}_K\] 
and \[\Sigma^2_{p_0}(W)= W \cap \Sec^2_{p_0}(W) = \{p \defin (0:0:1:0:\cdots:0)\}.\]
Moreover, $l(W\cap\Sec^2_{p}(W)) = 2$. So, $\Sec^2_{p}(W)$ is a line which intersects $W$ in one point $w$ with multiplicity $2$. It holds \[\langle w,\LL\rangle \cap W = \{w\} \mbox{ and } l(\langle w,\LL\rangle\cap W) = 3,\] meaning that $w$ `lies with length $3$ over its image $\pi(w)  \in Y$'. 
We now consider the PEIs corresponding to this projection. For this, we decompose $\pi_{\LL} = \tilde{\pi}\circ\pi'$, where \[\pi':\PP^{10}\bs\{p'\defin(0:0:0:0:1:0:\cdots:0)\} \lra \PP^9_K\] is given by $\tilde{S}' \defin K[x_0,x_1,x_2,x_3,x_5,\ldots,x_{10}] \inkl R$, and \[\tilde{\pi}:\PP^9_K\bs\{\tilde{p} \defin \pi'(p) = (0:0:0:1:0:\cdots:0)\} \lra \PP^8_K \] given by $S \inkl \tilde{S}'.$  Then, the PEIs of $I_W$ with respect to $p'$ are \[\K_0^{p'}(I_W) = I_W\cap \tilde{S}', \K_1^{p'}(I_W) = (\tilde{S}')_+, \K_2^{p'}(I_W) = \tilde{S}',\] where $I_W\cap \tilde{S}'$ is generated by $18$ quadrics in $S$ and $9$ quadrics and $1$ cubic in $\tilde{S}'\bs S$. $\K_1^{p_1}(I_W) = (\tilde{S}')_+$ means that $\pi'$ is an isomorphism in accordance with $p' \notin \Sec(W)$. For the PEIs of $\K_0^{p'}(I_W)$ with respect to $\tilde{p}$ we get \[\begin{array}{ll} \K^{\tilde{p}}_0(\K^{p'}_0(I_W)) = I_Y,\\ \K^{\tilde{p}}_1(\K^{p'}_0(I_W)) = (x_0,x_1^3,x_5,\ldots,x_{10}),\\ \K^{\tilde{p}}_2(K^{p'}_0(I_W)) = (x_0,x_1,x_5,\ldots,x_{10}),\\ \K^{\tilde{p}}_3(\K^{p'}_0(I_W)) = S.\end{array}\]
Looking at these PEIs, Corollary \ref{SemiIsoProj} tells us that $\pi_{\LL}(w) = (0:0:1:0:\cdots:0)$ is indeed the only point $q$ of $Y$ such that the length of the fibre $(\pi_{\LL}\res_W)^{-1}(q)$ is $3$; for every other point $q \in Y$ the length of the fibre is $l(W\cap\langle q,\LL\rangle) = 1$.\\
Finally, for example \[x_2x_5-x_3x_4 \in I_W\bs((I_W\cap K[p_1])R+(I_W \cap K[p'_1])R)\sat,\] so $p,p'$ is not a clever decomposition of $\LL$.\\
$\E \defin J_2S:_SQ = (x_5,\ldots,x_{10}) \ss S$ defines a projective plane $\EE = \PP^2 \ss \PP^8$. The intersection $Y\cap \EE$ is the quartic defined in $\EE$ by $\ol{Q} = x_1^3x_2$. Let $\LL_1$ be the projective line in $\EE$ defined by $x_1$, in $\PP^8$ by $\L_1 = (x_1, x_5, \ldots x_{10}) \ss S$. In $\PP^{10}$, $\L_1$ defines the projective three-space $\langle \LL,\LL_1\rangle$, and $I_W+\L_1R = (x_0x_4,x_2x_4-x_3^2,x_3x_4,x_4^2,x_1,x_5,\ldots,x_{10}) \ss R$, hence as a set $(\pi_{\LL}\res_W)^{-1}(\LL_1) = \langle \LL,\LL_1\rangle = V_{\PP^{10}}(x_1,x_3,x_4, \ldots, x_{10})$ is the ruling line $\widetilde{\LL}_1$ on $W$ which contains $w$. Moreover, $e_0(R/(I_W+\L_1R)) = 3$, that is '$\widetilde{\LL_1}$ lies with length $3$ over $\LL_1$'.
\end{ex}

\begin{ex}\label{ClDecEx}
Let $W = S(1,8) \ss \PP^{10}$ be as in Example \ref{7.4(E)}, but now consider the line $\LL = \PP^1 \ss \PP^{10}$ given by the ideal $\L = (x_0,x_1,x_2,x_4,\ldots,x_8,x_{10})$. Let $p,p' \in \LL$ be the closed points given by the ideals $p \defin (x_0,x_1,x_2,x_4,\ldots,x_{10})$, $p' \defin (x_0,\ldots,x_8,,x_{10}) \in \mProj(R)\cap\Var(\L)$, and let $\tilde{S} \defin K[p_1], \tilde{S}' \defin K[p'_1]$. Then a short computation using {\sc Singular} shows \[\left((I_W\cap \tilde{S})R + (I_W\cap \tilde{S}')R\right)\sat = I_W\sat = I_W,\] so $p,p'$ is a clever decomposition of $\LL$ with respect to $W$. Another consultation of {\sc Singular} gives \[\begin{array}{ll} 
\K_0^{\L\cap \tilde{S}}(I_W\cap \tilde{S}) = I_W\cap S = \K_0^{\L\cap \tilde{S}'}(I_W\cap \tilde{S}'), \\
\K_1^{\L\cap \tilde{S}}(I_W\cap \tilde{S}) = (x_0,x_1,x_2,x_4,\ldots,x_8), & \K_2^{\L\cap \tilde{S}}(I_W\cap \tilde{S}) = S,\\
\K_1^{\L\cap \tilde{S}'}(I_W\cap \tilde{S}') = (x_0,x_1,x_4,\ldots,x_8,x_{10}), & \K_2^{\L\cap \tilde{S}'}(I_W\cap \tilde{S}') = S.\end{array}\] 
Thus we can compute Table 2, giving us the ideals $\K_l^{\L\cap \tilde{S}}(I_W\cap \tilde{S}) + \K_k^{\L\cap \tilde{S}'}(I_W\cap \tilde{S}') \ss S$.
\begin{table}\centering\begin{tabular}{|c|ccc|}\hline
\backslashbox{l}{k} & 0 & 1 & 2\\ \hline
0 & $I_W\cap S$ & $(x_0,x_1,x_4,\ldots,x_8,x_{10})$ & $S$\\
1 & $(x_0,x_1,x_2,x_4,\ldots,x_8)$ & $S$ & $S$\\
2 & $S$ & $S$ & $S$\\ \hline
\end{tabular}\caption{$\K_l^{\L\cap \tilde{S}}(I_W\cap \tilde{S}) + \K_k^{\L\cap \tilde{S}'}(I_W\cap \tilde{S}') \ss S$}\end{table}
So we see that for $q = (0:0:1:0:\cdots:0), q' = (0:\cdots:0:1) \in Z$  it holds \[\K_0^{\L\cap \tilde{S}}(I_W\cap \tilde{S}) + \K_1^{\L\cap \tilde{S}'}(I_W\cap \tilde{S}') = q\] and \[\K_1^{\L\cap \tilde{S}}(I_W\cap \tilde{S}) + \K_0^{\L\cap \tilde{S}'}(I_W\cap \tilde{S}') = q',\] i.e. $\length((\pi_{\LL}\res_Z)^{-1}(q)) = \length((\pi_{\LL}\res_Z)^{-1}(q')) = 2$. Indeed, we can compute \[\Sigma^2_{p}(Z) = \{(\pi_{\LL}\res_Z)^{-1}(q)\} \mbox{ and } \Sigma^2_{p'}(Z) = \{(\pi_{\LL}\res_Z)^{-1}(q')\}.\]
\end{ex}

\end{document}